\documentclass[11pt]{article}
\usepackage{amsfonts,amssymb,amsmath}

\newcommand{\bP}{{\bf P}}

\newcommand{\bbZ}{\mathbb Z}
\newcommand{\bbR}{\mathbb R}
\newcommand{\bbN}{\mathbb N}

\newcommand{\bbE}{\mathbb E}

\newcommand{\bbP}{\mathbb P}

\usepackage[english]{babel}
\newtheorem{theo}{Theorem}[section]
\newtheorem{pr}{Proposition}[section]

\newtheorem{lem}{Lemma}[section]
\newtheorem{defn}{Definition}[section]

\newcommand{\Cov}{\mathop{\rm Cov}}

\newcommand{\epreuve}{\hspace{\fill}$\bigtriangleup$}

\newcommand{\NB}{{\par\noindent{\bf Remark~:~}}}

\begin{document}

\title{{\bf Limit theorem for random walk in weakly dependent random scenery}}

\author{Nadine {\bf Guillotin-Plantard}\footnote{Universit\'e de Lyon, Institut Camille Jordan,
B\^atiment Braconnier, 43 avenue du 11 novembre 1918,  69622
Villeurbanne Cedex, France. {\it E-mail:}
nadine.guillotin@univ-lyon1.fr} and Cl\'ementine {\bf Prieur}\footnote{INSA Toulouse, Institut Math\'ematique de Toulouse, Equipe de Statistique et Probabilit\'es,
135 avenue de Rangueil, 31077 Toulouse Cedex 4, France.
{\it E-mail:} Clementine.Prieur@insa-toulouse.fr}}

\maketitle

\begin{abstract}
Let $S=(S_k)_{k\geq 0}$ be a random walk on $\bbZ$ and $\xi=(\xi_{i})_{i\in\bbZ}$ a stationary random sequence of centered random variables, independent of $S$. We consider a random walk in random scenery that is the sequence of random variables $(\Sigma_n)_{n\geq 0}$ where
$$\Sigma_n=\sum_{k=0}^n \xi_{S_k},\ n\in\bbN.$$ 
Under a weak dependence assumption on the scenery $\xi$ we prove a functional limit theorem generalizing Kesten and Spitzer's theorem (1979).
\end{abstract}

\medskip


\vspace{7cm}



\newpage
\section{Introduction}

Let $X=(X_{i})_{i\geq 1}$ be a sequence of independent and identically distributed random vectors with values in $\bbZ^{d}$.
We write $$S_{0}=0,\ \ S_{n}=\sum_{i=1}^{n} X_{i}\ \mbox{for}\  n\geq 1$$
for the $\bbZ^{d}$-random walk $S=(S_n)_{n\in\bbN}$ generated by the family $X$.
Let $\xi=(\xi_{x})_{x\in \bbZ^{d}}$ be a family of real random variables, independent of $S$. The sequence $\xi$ plays the role of the random scenery. The random walk in random scenery (RWRS) is the process defined by 
$$\Sigma_n=\sum_{k=0}^n \xi_{S_k},\ n\in\bbN. $$ 
RWRS was first introduced in dimension one by Kesten and Spitzer (1979) and Borodin (1979) in order to construct new self-similar stochastic processes. Functional limit theorems for RWRS were obtained under the assumption that the random variables $\xi_x, x\in\bbZ^{d}$ are independent and identically distributed. 
For $d=1$, Kesten and Spitzer (1979) proved that when $X$ and
$\xi$ belong to the domains of attraction of different stable laws
of indices $1<\alpha\leq 2$ and $0<\beta\leq 2$, respectively,
then there exists $\delta>\frac{1}{2}$ such that
$\big(n^{-\delta}\Sigma_{[nt]}\big)_{t\geq 0}$ converges weakly as $n\rightarrow
\infty$ to a self-similar process with stationary increments,
$\delta$ being related to $\alpha$ and $\beta$
 by $\delta=1-\alpha^{-1}+(\alpha\beta)^{-1}$.
The case $0<\alpha<1$ and $\beta$ arbitrary is easier; they showed
then that $\big(n^{-\frac{1}{\beta}}\Sigma_{[nt]}\big)_{t\geq 0}$ converges
weakly, as $n\rightarrow \infty$, to a stable process with index
$\beta$. Bolthausen (1989) gave a method to solve the case $\alpha=1$ and $\beta=2$ and especially, he proved
that when $(S_{n})_{n\in \bbN}$ is a recurrent $\bbZ^{2}$-random
walk, $\big((n\log n)^{-\frac{1}{2}}\Sigma_{[nt]}\big)_{t\geq 0}$ satisfies a functional
central limit theorem. For an arbitrary transient $\bbZ^{d}$-random walk,
$n^{-\frac{1}{2}}\Sigma_{n}$ is asymptotically normal (see \cite{Spi} p. 53). 
Maejima (1996) generalized the result of Kesten and
Spitzer (1979) in the case where $(\xi_x)_{x \in \bbZ}$ are
i.i.d. $\bbR^d$-valued random variables which belong to the domain
of attraction of an operator stable random vector with exponent
$B$. If we denote by $D$  the linear operator on $\bbR^d$ defined by $D=(1-\frac{1}{\alpha})I +\frac{1}{\alpha}B$, he proved that  $\big(n^{-D}\Sigma_{[nt]}\big)_{t\geq 0}$ converges weakly to an operator self similar  with exponent $D$ and having stationary increments.

One-dimensional random walks in random scenery recently arose in the study of random walks evolving on oriented versions of $\bbZ^{2}$ (see Guillotin-Plantard and Le Ny (2007, 2008)) as well as in the context of charged polymers (see Chen and Khoshnevisan (2008)). The understanding of these models in the case where the orientations or the charges are not independently distributed requires functional limit theorems for $\bbZ$-random walk in correlated random scenery. To our knowledge, only the case of strongly correlated stationary random sceneries 
has been studied by Lang and Xanh (1983). 
In their paper, the increments of the random walk $S$ are assumed to belong to the domain of attraction of a
non-degenerate stable law of index $\alpha, 0 < \alpha\leq 2$. They further suppose that the
scenery $\xi$ satisfies the non-central limit theorem of Dobrushin and Major (1979) with
a scaling factor $n^{-d+(\beta k)/2}$, $\beta k<d$. Under the assumption $\beta k<\alpha$, it is proved
that $\big(n^{-1+\beta k/ (2\alpha)}\Sigma_{[nt]}\big)_{t\geq 0}$ converges weakly as $n\rightarrow +\infty$ to a self-similar process with stationary increments, which can be represented as a multiple Wiener-It\^{o} integral of a random function.
Our aim is to study the intermediary case of a stationary random scenery $\xi$ which satisfies a weak dependence condition introduced in Dedecker {\it et al.} (2007) and to prove Kesten and Spitzer's theorem under this new assumption. In Guillotin-Plantard \& Prieur (2008) the case of a transient $\bbZ$-random walk was considered and a central limit theorem for the sequence $(\Sigma_n)_{n\in\bbN}$ was proved. In this paper the one-dimensional random walk will be assumed to be recurrent. 
  
Our paper is organized as follows: In Section \ref{def}, we introduce the dependence setting under
which we work in the sequel. In Section \ref{res} we introduce in details our model 
and give the main result. In Section \ref{time} properties of the local time of the random walk 
are given as well as the ones of the intersection local time.
Models for which we can compute bounds
for our dependence coefficients are presented in Section
\ref{examples}. Finally, the proof of our theorem is given in the last section.

\section{Weak dependence conditions}\label{def}
In this section, we recall the definition of the dependence coefficients
which we will use in the sequel. They have first been introduced in Dedecker
{\it et al.} (2007). Our weak dependence condition will be less restrictive than the mixing one. The reader interested in this question would find more details in Guillotin-Plantard \& Prieur (2008).

On the Euclidean space $\bbR^m$, we define the metric
$$d_1(x,y)= \sum_{i=1}^m |x_i - y_i|.$$
Let $\Lambda=\bigcup_{m\in\bbN^{*}}\Lambda_m$ where $\Lambda_m$ is the set of Lipschitz functions $f:\bbR^m\rightarrow \bbR$ with respect to the metric $d_1$.
If $f\in\Lambda_m$, we denote by $\mbox{\rm Lip}(f):=\sup_{x\neq y} \frac{|f(x)-f(y)|}{d_1(x,y)}$ the Lipschitz modulus of $f$. The set of functions $f\in \Lambda$ such that $\mbox{\rm Lip}(f)\leq 1$ is denoted by $\tilde\Lambda$.\\*

\begin{defn}\label{coeff}
 Let $\xi$ be a $\bbR^m$-valued random variable defined on a probability space
$(\Omega,{\cal A},\bbP)$, assumed to be square integrable.
For any $\sigma$-algebra ${\cal M}$ of ${\cal A}$, we define the $\theta_2$-dependence coefficient
\begin{equation}\label{tet2}
\theta_2({\cal M},\xi)=\sup \{\|\bbE(f (\xi)|{\cal M})-\bbE(f (\xi))\|_2 \, , \;
f\in \tilde\Lambda\} \, .
\end{equation}
\end{defn}

We now define the coefficient $\theta_{k,2}$ for a sequence of
$\sigma$-algebras and a sequence of $\mathbb{R}$-valued random variables.

\begin{defn}\label{coeffsuite}
Let $(\xi_i)_{i\in \bbZ}$ be a sequence of square integrable random variables
valued in $\mathbb{R}$. Let $({\cal M}_i)_{i\in\bbZ}$ be a sequence of
$\sigma$-algebras of ${\cal A}$. For any $k \in \mathbb{N}^*\cup \{\infty\}$
and $n \in \mathbb{N}$, we define
$$\theta_{k,2}(n)= \max_{1\leq l\leq k}
\frac 1 l \sup\{\theta_{2}({\cal M}_p, (\xi_{j_1}, \ldots ,\xi_{j_l})), p+n \leq j_1< \ldots <j_l\}$$ and
$$\theta_2(n)= \theta_{\infty,2} (n)= \sup_{k\in \mathbb{N}^*} \theta_{k,2}(n) \, .$$
\end{defn}
\begin{defn}
Let $(\xi_i)_{i\in \bbZ}$ be a sequence of square integrable random variables
valued in $\mathbb{R}$. Let $({\cal M}_i)_{i\in\bbZ}$ be a sequence of
$\sigma$-algebras of ${\cal A}$. The sequence $(\xi_i)_{i \in \bbZ}$ is said
to be $\theta_2$-weakly dependent with respect to $({\cal M}_i)_{i\in\bbZ}$
if $\theta_2(n) \xrightarrow[n \rightarrow + \infty]{} 0$.
\end{defn}

\NB
Replacing the
$\| \cdot \|_2$ norm in (\ref{tet2}) by the $\| \cdot \|_1$ norm, we get the
$\theta_1$ dependence coefficient first introduced by Doukhan \& Louhichi
(1999).

\section{Model and results}\label{res}
Let $S=(S_k)_{k\geq 0}$ be a $\bbZ$-random walk ($S_0=0$) whose increments $(X_i)_{i\geq 1}$ are centered and square integrable. We denote by $\bP_{X_1}$ the law of the random variable $X_1$. For any $q\in \bbN^{*}$ such that $\bP(X_1\in [-q,q])$ is non zero, we define the probability measure on $\bbZ$
$${\bP}_{q}=\frac{{\bP_{X_1}}\big|_{[-q,q]}}{\bP(X_1\in [-q,q])}.$$
The random walk $S$ is said to satisfy the property $({\bf P})$ if there exists $q\in \bbN^{*}$ such that
$$\{x\in\bbZ;\ \  \exists n,\ {\bP}_{q}^{(\star n)} (x) >0\}=\bbZ.$$ 
In particular, if there exists some $q\in \bbN^{*}$ such that the random walk associated to ${\bP}_{q}$ is centered and aperiodic then $S$ satisfies the property $({\bf P})$. For instance, the simple random walk on $\bbZ$ verifies ${\bf (P)}$. Let $\xi=(\xi_i)_{i\in\mathbb{Z}}$ be a sequence of centered real random variables. The sequences $S$ and $\xi$ are defined on a same probability space denoted by $(\Omega, {\cal F},\bP)$ and are assumed to be independent.
We are interested in the asymptotic behaviour of the following sum
$$\Sigma_n=\sum_{k=0}^n \xi_{S_k} \, .$$
The case where the $\xi_i$'s are independent and identically distributed random variables with positive variance has been considered by Kesten \& Spitzer (1979) and Borodin (1979). Consider a standard Brownian motion $(B_{t})_{t\ge 0}$,
denote by $(L_{t}(x))_{t \ge 0}$ its corresponding local time at
$x\in\bbR$ and introduce a pair of independent Brownian motions
$(Z_{+}(x), Z_{-}(x)), x\geq 0$ defined on the same probability
space as $(B_{t})_{t\ge 0}$ and independent of him. The following
process is well-defined for all $t \geq 0$:
\begin{equation}\label{th}
\Delta_{t}=\int_{0}^{\infty}L_{t}(x)dZ_{+}(x)+\int_{0}^{\infty}L_{t}(-x)dZ_{-}(x).
\end{equation}
It was proved by Kesten and Spitzer (1979) that this
process has a self-similar continuous version  of index
$\frac{3}{4}$, with stationary increments. We denote
$\stackrel{\mathcal{D}}{\Longrightarrow}$ for a  convergence in
the space of c\`adl\`ag functions $\mathcal{D}([0,\infty),\bbR)$
endowed with the Skorohod topology.
\begin{theo}\label{thm11} \mbox{\bf [Kesten and Spitzer (1979)]}
Assume that the $\xi_i$'s are independent and identically distributed with positive variance $\sigma^2>0$. Then,
 \begin{equation}
  \Big(\frac{1}{n^{3/4}} {\Sigma}_{[nt]} \Big)_{t \geq 0}
\; \stackrel{\mathcal{D}}{\Longrightarrow} (\sigma \Delta_t)_{t \geq 0}.
\end{equation}
\end{theo}
A simple proof of this theorem was proposed by Cadre (1995) (Section 2.5.a.) using a weak limit theorem for stochastic integrals (Theorem 1.1. in Kurtz \& Protter (1991)). 
Applying Cadre's method, Theorem \ref{thm11} can be extended to any scenery given by a stationary and ergodic sequence of square integrable martingale differences. Then, a natural idea is to generalize the result to any stationary and ergodic sequence $\xi$ of square integrable random variables as it was done for the central limit theorem. Under suitable assumptions on the sequence, for instance 
the convergence of the series $\sum_{k=0}^{\infty} \mathbb{E}(\xi_k|{\cal M}_0)$ in ${\mathbb{L}}^2$, the scenery $\xi$ is equal to a martingale differences sequence modulo a coboundary term and satisfies a Donsker theorem. However, the RWRS associated to the coboundary term (if it is non zero) is not negligible. It can be proved that the ${\mathbb{L}}^2-$norm of this sum correctly normalized by $n^{3/4}$ converges to some positive constant. 

In order to weaken the assumptions on the field $\xi$ we introduce a sequence $({\cal M}_i)_{i \in \mathbb{Z}}$ of
$\sigma$-fields of ${\cal F}$ defined by
$$\mathcal{M}_i=\sigma\left(\xi_j \, , \; j \leq i\right)\, , \;
i\in \mathbb{Z} \, .$$ In the sequel, the dependence coefficients will be
defined with respect to the sequence of $\sigma$-fields $({\cal
M}_i)_{i \in \mathbb{Z}}$. 
\begin{theo}\label{thmix0}
\hspace{1em}\\
Assume that the following conditions are satisfied~:
\begin{itemize}
\item[$(A_0)$] The random walk $S$ satisfies the property ${\bf (P)}$.
\item[$(A_1)$]
$\xi=\{\xi_i\}_{i \in \mathbb{Z}}$ is a stationary sequence of square integrable random variables.
\item[$(A_2)$] $\theta_2^{\xi}(\cdot)$ is bounded above by a
  non-negative function $g(\cdot)$ such that
\begin{itemize}
\item[] $x \mapsto x^{3/2} \ g(x)$ is non-increasing,
\item[] $\exists \ 0 < \varepsilon < 1 \ , \, \displaystyle\sum_{i=0}^{\infty}
  2^{\frac{3i}{2}}g(2^{i\varepsilon}) < \infty $.
\end{itemize}
\end{itemize}
Then, as $n$ tends to infinity, 
\begin{equation}
  \Big(\frac{1}{n^{3/4}} {\Sigma}_{[nt]} \Big)_{t \geq 0}
\; \stackrel{\mathcal{D}}{\Longrightarrow} \displaystyle\sqrt{\sum_{i\in\bbZ} \mathbb{E}(\xi_0 \xi_i)}\  (\Delta_t)_{t \geq 0}.
\end{equation}
\end{theo}

\NB
Assumptions $(A_1)$ and $(A_2)$ imply that 
\begin{equation}\label{imp}
\forall \lambda\in [0,1 / 2[,\ \  \sum_{k\in\bbZ} |k|^{\lambda} |\mathbb{E}(\xi_0\xi_k)|<+\infty
\end{equation}
Indeed, this sum is equal to
$$ \mathbb{E}(\xi_0^2) + 2 \sum_{k=1}^{\infty} |k|^{\lambda} |\mathbb{E}(\xi_0 \xi_k)|$$
and for any $k\geq 1$, from Cauchy-Schwarz inequality, we get 
\begin{eqnarray*}
|\mathbb{E}(\xi_0 \xi_k)| &=& | \mathbb{E}\big(\xi_0 \mathbb{E}(\xi_k|{\cal M}_0)\big)| \\
&\leq & ||\xi_0||_2\ \theta_2^{\xi}(k)\\
&\leq & ||\xi_0||_2\ g(k)
\end{eqnarray*}
The result $(\ref{imp})$ follows by remarking that
 $$\sum_{k=1}^{\infty} |k|^{\lambda} g(k)\leq g(1) \sum_{k=1}^{\infty} \frac{1}{k^{3/2-\lambda}}$$
which is finite for every $\lambda\in [0,1 / 2[$.

\NB
 If $\theta_2^{\xi}(n)=\mathcal{O}\left(n^{-a}\right)$ for some
  positive $a$, condition $(A_3)$ holds for $a>3/2$.

\section{Properties of the occupation times of the random walk}\label{time}
The random walk $S=(S_k)_{k\geq 0}$ is defined as in the previous section and is assumed to verify the property ${\bf (P)}$.
The local time of the random walk is defined for every $i\in\bbZ$ by
$$N_{n}(i)= \sum_{k=0}^{n} {\bf 1}_{\{S_{k}=i\}} \; .$$
The local time of self-intersection at point $i$ of the random walk $(S_n)_{n\geq 0}$ is defined  by
$$\alpha(n,i)= \sum_{k,l=0}^{n} {\bf 1}_{\{S_{k}-S_{l}=i\}}.$$
The stochastic properties of the sequences $(N_n(i))_{n\in\bbN,i\in\bbZ}$ and $(\alpha(n,i))_{n\in\bbN,i\in\bbZ}$ are well-known when the random walk $S$ is strongly aperiodic. A random walk who satisfies the property ${\bf (P)}$ is not strongly aperiodic in general. However, a local limit theorem for the random walks satisfing ${\bf (P)}$ was proved by Cadre (2005) (see Lemma 2.4.5., p. 70), then it is not difficult to adapt the proofs of the strongly aperiodic case to our setting: for assertion (i) see Lemma 4 in Kesten and Spitzer (1979), for (ii)-(a) see Lemma 3.1 in Dombry and Guillotin-Plantard (2008). Result (ii)-(b) is obtained from Lemma 6 in Kesten and Spitzer (1979); details are omitted. Assertion (iii) is an adaptation to dimension one of Lemma 2.3.2 of Cadre's thesis (1995). 
\begin{pr}\label{cadre}
\begin{itemize}
\item[(i)] The sequence $ n^{-3/4} \displaystyle\max_{i \in \mathbb{Z}} N_n(i)$ converges in probability to 0.
\item[(ii)] 
\begin{itemize}
\item[(a)] For any $p\in [1,+\infty)$, there exists some constant $C$ such that for all $n\geq 1$, 
$$\bbE\left(\alpha(n,0)^p\right)\leq C n^{3p/2}.$$
\item[(b)] For any $m\geq 1$, for any real $\theta_1,\ldots,\theta_m$, for any $0\leq t_1\leq  \ldots\leq t_m$, the sequence 
$$n^{-3/2} \sum_{i\in\bbZ}\Big(\sum_{k=1}^m \theta_k N_{[nt_k]}(i)\Big)^2$$ 
converges in distribution to 
$$\int_{\mathbb{R}} \Big(\sum_{k=1}^{m} \theta_k L_{t_k}(x)\Big)^2 \, dx$$  
where $(L_t(x))_{t\geq 0;x\in\mathbb{R}}$ is the local time of the real Brownian motion $(B_t)_{t\geq 0}$. 
\end{itemize}
\item[(iii)] For every $\lambda\in\ ]0,1[$, there exists a constant $C$ such that for any $i,j\in \mathbb{Z}$, 
$$ \| \alpha(n,i) - \alpha(n,j) \|_2\leq C n^{(3-\lambda)/2} |i-j|^{\lambda} . $$
\end{itemize}
\end{pr}

\section{Examples}\label{examples}
In this section, we present examples for which we can compute upper bounds
for $\theta_{2}(n)$ for any $n \geq 1$. We refer to Chapter 3 in Dedecker {\it
et al} (2007) and references therein for more details.
\subsection{Example 1: causal functions of stationary sequences}\label{ex1}
Let $(E,\mathcal{E},\mathbb{Q})$ be a probability space.
Let $(\varepsilon_i)_{i \in \mathbb{Z}}$ be a stationary sequence of random
variables with values in a measurable space $\mathcal{S}$. Assume that there
exists a real valued function $H$ defined on a subset of
$\mathcal{S}^{\mathbb{N}}$, such that $H(\varepsilon_0,\varepsilon_{-1},
\varepsilon_{-2}, \ldots,)$ is defined almost surely. The stationary
sequence $(\xi_n)_{n \in \mathbb{Z}}$ defined by $\xi_n=H(\varepsilon_n,
\varepsilon_{n-1}, \varepsilon_{n-2}, \ldots)$ is called a causal function
of $(\varepsilon_i)_{i \in \mathbb{Z}}$.

Assume that there exists a stationary sequence $({\varepsilon_i}')_{i \in
  \mathbb{Z}}$ distributed as $(\varepsilon_i)_{i \in \mathbb{Z}}$ and
  independent of $(\varepsilon_i)_{i \leq 0}$. Define $\xi_n^{*}=H({\varepsilon_n}',
{\varepsilon_{n-1}}', {\varepsilon_{n-2}}', \ldots)$. Clearly, $\xi_n^{*}$ is
  independent of $\mathcal{M}_0=\sigma(\xi_i \, , \; i \leq 0)$ and
  distributed as $\xi_n$. Let $(\delta_2(i))_{i >0}$ be a non increasing
  sequence such that
\begin{equation}\label{delta2}
\left\| \mathbb{E}\left(\left|\xi_i-\xi_i^{*}\right| \, | \,
  \mathcal{M}_0\right)\right\|_2 \leq \delta_2(i) \, .
\end{equation}
Then the
  coefficient $\theta_2$ of the sequence $(\xi_n)_{n \geq 0}$ satisfies
\begin{equation}\label{bound}
\theta_2(i) \leq \delta_2(i) \, .
\end{equation}
Let us consider the particular case where the sequence of innovations $(\varepsilon_i)_{i \in
  \mathbb{Z}}$ is absolutely regular in the sense of Volkonskii \& Rozanov
(1959). Then, according to Theorem 4.4.7 in Berbee (1979), if $E$ is
rich enough, there exists $(\varepsilon_i')_{i \in \mathbb{Z}}$ distributed
as $(\varepsilon_i)_{i \in \mathbb{Z}}$ and independent of
$(\varepsilon_i)_{i \leq 0}$ such that $$\mathbb{Q}(\varepsilon_i \neq
\varepsilon_i' \textrm { for some } i \geq k \, | \, \mathcal{F}_0)=\frac12
\left\| \mathbb{Q}_{\tilde{\varepsilon}_k|\mathcal{F}_0} -
\mathbb{Q}_{\tilde{\varepsilon}_k}\right\|_v ,$$
where $\tilde{\varepsilon}_k=(\varepsilon_k, \varepsilon_{k+1}, \ldots)$,
$\mathcal{F}_0=\sigma(\varepsilon_i \, , \, i \leq 0)$, and $\| \cdot \|_{v}$ is
  the variation norm. In particular if the
sequence $(\varepsilon_i)_{i \in \mathbb{Z}}$ is idependent and identically
distributed, it suffices to take $\varepsilon_i'=\varepsilon_i$ for $i>0$
and $\varepsilon_i'-\varepsilon_i''$ for $i \leq 0$, where
$(\varepsilon_i'')_{i \in \mathbb{Z}}$ is an independent copy of
$(\varepsilon_i)_{i \in \mathbb{Z}}$.

\noindent {\bf Application to causal linear processes:}\\
In that case, $\xi_n=\sum_{j \geq 0} a_j \varepsilon_{n-j}$, where $(a_j)_{j
  \geq 0}$ is a sequence of real numbers. We can choose
$$\delta_2(i) \geq \|\varepsilon_0-\varepsilon_0'\|_2 \sum_{j \geq i} |a_j| + \sum_{j=0}^{i-1} |a_j| \|\varepsilon_{i-j}- \varepsilon_{i-j}'\|_2 \, .
$$
From Proposition 2.3  in Merlev\`ede \& Peligrad (2002), we obtain
that
$$
  \delta_2(i) \leq  \|\varepsilon_0-\varepsilon_0'\|_2 \sum_{j \geq i} |a_j| + \sum_{j=0}^{i-1} |a_j| \,
\Big (2^2 \int_0^{\beta(\sigma(\varepsilon_k, k \leq 0), \sigma(\varepsilon_k, k \geq i-j))} Q^2_{\varepsilon_0}(u)\Big)^{1/2} du
 ,
$$
where $Q_{\varepsilon_0}$ is the generalized inverse of the tail function
$x \mapsto \mathbb{Q}(|\varepsilon_0|>x)$.

\subsection{Example 2: iterated random functions}
Let $(\xi_n)_{n \geq 0}$ be a real valued stationary Markov chain,
such that
$
\xi_n=F(\xi_{n-1}, \varepsilon_n)
$
for some measurable function $F$ and some independent and identically
distributed sequence $(\varepsilon_i)_{i > 0}$ independent of $\xi_0$.
Let $\xi^*_0$ be a random variable distributed as $\xi_0$ and independent
of
$(\xi_0, (\varepsilon_i)_{i >0})$.
Define
$
\xi^*_n=F(\xi^*_{n-1}, \varepsilon_n) \, .
$
The sequence $(\xi^*_n)_{n \geq 0}$ is distributed as $(\xi_n)_{n \geq  0}$
and independent of
$\xi_0$. Let ${\cal M}_i=\sigma(\xi_j , 0 \leq j \leq i)$.
As in Example 1,
define the sequence
$(\delta_2(i))_{i > 0}$ by (\ref{delta2}).
The coefficient  ${\theta}_2$ of
the sequence $(\xi_n)_{n \geq 0}$ satisfies the bound (\ref{bound}) of Example 1.

Let $\mu$ be the distribution of
$\xi_0$ and $(\xi_n^x)_{n \geq 0}$ be the chain starting from $\xi_0^x=x$.
With these notations, we can choose $\delta_2(i)$ such that
$$\delta_2(i) \geq \| \xi_i-\xi^*_i\|_2=\left(\int\int \|
|\xi_i^x-\xi_i^y\|_2^2 \mu(dx) \mu(dy) \right)^{1/2} \, .$$
For instance, if there exists a sequence
$(d_2(i))_{i \geq 0}$ of positive numbers
such that
$$
\|\xi_i^x-\xi_i^y\|_2 \leq d_2(i) |x-y|  ,
$$
then we can take
$ \delta_2(i) =d_2(i) \| \xi_0-\xi_0^* \|_2$.
For example, in the usual case where
$ \|F(x, \varepsilon_0)-F(y, \varepsilon_0)\|_2 \leq \kappa |x-y| $
for some $\kappa < 1$, we can
take $d_2(i)= \kappa^i$.

An important example is $\xi_n=f(\xi_{n-1}) + \varepsilon_n$
for some $\kappa$-Lipschitz function $f$. If $\xi_0$
has a moment of order $2$, then
$
\delta_2(i) \leq \kappa^i \|\xi_0-\xi_0^*\|_2 \, .
$

\subsection{Example 3: dynamical systems on $[0,1]$}
Let $I=[0,1]$, $T$ be a map from $I$ to $I$ and define $X_i=T^i$. If
$\mu$ is invariant by $T$, the sequence
$(X_i)_{i\geq0}$ of random variables from $(I,\mu)$ to $I$ is strictly stationary.

For any finite measure $\nu$ on $I$, we use the notations
$\nu(h)=\int_I h(x) \nu(dx)$.
For any finite signed measure $\nu$ on $I$, let $\|\nu\|=|\nu|(I)$ be the total variation
of $\nu$.
 Denote by $\|g\|_{1, \lambda}$ the ${\mathbb L}^1$-norm with respect to
the Lebesgue measure $\lambda$ on $I$. \medskip

\noindent{\bf Covariance inequalities.} In many interesting cases, one can prove that,
 for any $BV$ function $h$ and any $k$ in ${\mathbb L}^1(I, \mu)$,
\begin{equation}\label{sd1}
|\Cov(h(X_0),k(X_n))| \leq a_n  \|k(X_n)\|_1(\|h\|_{1, \lambda}+ \|dh\|) \, ,
\end{equation}
for some nonincreasing sequence $a_n$ tending to zero as $n$ tends to
infinity.

\medskip

\noindent{\bf Spectral gap.}
Define the operator $\mathcal{L}$ from ${\mathbb L}^1(I, \lambda)$ to ${\mathbb L}^1(I, \lambda)$ $via$ the equality
$$
\int_0^1 \mathcal{L}(h)(x) k(x) d \lambda (x) = \int_0^1h(x) (k \circ T)(x) d \lambda (x) \ \, \text{ where $h \in {\mathbb L}^1(I, \lambda)$ and
$k \in {\mathbb L}^\infty(I, \lambda)$}.
$$
The operator $\mathcal{L}$ is called the Perron-Frobenius operator of $T$.
In many interesting cases, the spectral analysis of ${\cal L}$  in the Banach space of $BV$-functions
equiped with the norm $\|h\|_v=\|dh\|+\|h\|_{1, \lambda}$ can be done by using the Theorem of
Ionescu-Tulcea and Marinescu (see Lasota and Yorke (1974) and Hofbauer and Keller (1982)).
Assume that $1$ is a simple eigenvalue of ${\cal L}$ and that the rest
of the spectrum is contained in a closed disk of radius strictly smaller than one.
Then there exists
a unique $T$-invariant absolutely continuous probability $\mu$ whose density $f_\mu$ is $BV$, and
\begin{equation}\label{sg1}
 {\cal L}^n(h)=\lambda(h) f_\mu + \Psi^n(h) \quad \text{with} \quad
\|\Psi^n(h)\|_v \leq K\rho^n \|h\|_v.
\end{equation}
for some $0\leq \rho <1$ and $K>0$.
Assume moreover that:
\begin{equation}\label{sg2}
\text{$I_*=\{f_\mu \neq 0\}$ is an interval, and there exists $\gamma>0$ such that $f_\mu > \gamma^{-1}$ on $I_*$.}
\end{equation}
Without loss of generality assume that $I_*=I$ (otherwise, take the restriction
to $I_*$ in what follows).
Define now the Markov kernel associated to $T$ by
\begin{equation}\label{Mk}
P(h)(x)=\frac{\mathcal{L}(f_\mu h)(x)}{f_\mu(x)} .
\end{equation}
It is easy to check (see for instance Barbour {\it et al.} (2000)) that $(X_0,X_1, \ldots , X_n)$
has the same distribution as $(Y_n,Y_{n-1},\ldots,Y_0)$ where $(Y_i)_{i \geq 0}$
is a stationary Markov chain with invariant distribution $\mu$ and transition kernel $P$.
Since $\|fg\|_\infty \leq \|fg\|_v \leq 2\|f\|_v \|g\|_v$, we infer that, taking $C=2K\gamma (\|df_\mu\|+1)$,
\begin{equation}\label{sg3}
 P^n(h)=\mu(h) + g_n \quad \text{with} \quad
\|g_n\|_\infty \leq C\rho^n \|h\|_v.
\end{equation}
This estimate implies   (\ref{sd1}) with $a_n=C\rho^n$ (see Dedecker \&
Prieur, 2005).

\medskip

\noindent {\bf Expanding maps:} Let  $([a_i, a_{i+1}[)_{1\leq i \leq N}$ be a finite partition of $[0,1[$.
We make the same assumptions on $T$ as in Collet {\it et al} (2002).
\begin{enumerate}
\item For each $1 \leq j \leq N$, the restriction $T_j$ of $T$ to $]a_j, a_{j+1}[$ is strictly monotonic
and can be extented to a function $\overline T_j$ belonging to $C^2([a_j, a_{j+1}])$.
\item Let $I_n$ be the set where $(T^n)'$ is defined.
There exists  $A>0$ and $s>1$ such that $ \inf _{x \in I_n} |(T^n)'(x)|>As^n$.
\item The map $T$ is topologically mixing: for any two nonempty open sets $U,V$, there exists $n_0\geq 1$
such that $T^{-n}(U)\cap V \neq \emptyset$ for all $n \geq n_0$.
\end{enumerate}
If $T$ satisfies 1., 2. and 3., then (\ref{sg1}) holds. Assume furthermore
that (\ref{sg2}) holds (see Morita
(1994) for sufficient conditions). Then, arguing as in Example 4 in Section
7 of Dedecker \& Prieur (2005), we can prove that for the Markov chain
$(Y_i)_{i \geq 0}$ and the $\sigma$-algebras $\mathcal{M}_i=\sigma(Y_j \, ,
\, j \leq i)$, there exists a positive constant $C$ such that
$\theta_2(i) \leq C \rho^i $.

\medskip

\section{Proof of Theorem \ref{thmix0}}
The proof of Theorem \ref{thmix0} is decomposed in two parts: first, we prove the convergence of the finite-dimensional distributions of the process $(n^{-3/2} \Sigma_{[nt]})_{t\geq 0}$, then its tightness in the Skorohod space ${\cal D}([0,+\infty[)$.
\\*
\\*
\noindent{\bf Proof of the convergence of the finite-dimensional distributions:}\\* 
Since the random variable $\Sigma_n$ can be rewritten as the sum
$$\sum_{i\in\bbZ} N_n(i) \xi_{i}$$ 
where $N_n(i)$ is the local time of the random walk $S$ at point $i$, it is enough to prove that for every $m\geq 1$, for any real $\theta_1,\ldots,\theta_m$, for any $0\leq t_1\leq t_2\leq \ldots\leq t_m$, the sequence
$$\frac{1}{n^{3/2}} \sum_{k=1}^m \theta_k \Sigma_{[nt_k]}=\frac{1}{n^{3/2}} \sum_{i\in\bbZ} \Big(\sum_{k=1}^m \theta_k N_{[nt_k]}(i)\Big) \xi_i $$
converges in distribution to the random variable
$$\sqrt{\sum_{i\in\bbZ} \mathbb{E}(\xi_0\xi_i)}\ \sum_{k=1}^m \theta_k \Delta_{t_k}.$$ 
We only prove the convergence of one-dimensional distributions. The general case is obtained by replacing $N_n(i)$
by the linear combination $\sum \theta_k N_{[nt_k]}(i)$ in the computations.\\*
Let ${\cal G}=\sigma(S_k, k\geq 0)$ be the $\sigma-$field generated by the random walk $S$. For any $n\in\bbN$ and any $i\in\bbZ$, we denote by $X_{n,i}$ the random variable $N_n(i) \xi_i$.\\*
We first use a classical truncation argument. For any $M>0$, we define~:
\begin{equation*}
\varphi_{M}~:
\begin{cases}
\mathbb{R} \rightarrow \mathbb{R}\\
x \mapsto \varphi_{M}(x)=(x \wedge M) \vee (-M)
\end{cases}
\end{equation*}
and
\begin{equation*}
\varphi^M~:
\begin{cases}
\mathbb{R} \rightarrow \mathbb{R}\\
x \mapsto \varphi^M(x)=x-\varphi_{M}(x) \, .
\end{cases}
\end{equation*}
We now prove the following Lindeberg condition~:
\begin{equation}\label{lind}
n^{-3/2}   \sum_{i\in\mathbb{Z}} \mathbb{E} \left( \left( \varphi^{\varepsilon
  n^{3/4}}(X_{n,i})\right)^2 \right) \xrightarrow[n \rightarrow + \infty]{}
  0 \, .
\end{equation}
We have, for $\varepsilon>0$ fixed, for $n$ large enough,
$$\begin{array}{rcl}
& & n^{-3/2}\displaystyle\sum_{i\in\mathbb{Z}} \mathbb{E} \left( \left( \varphi^{\varepsilon
  n^{3/4}}(X_{n,i})\right)^2\Big|{\cal G}\right)\\   
  &\leq & \mathbb{E} \left( \xi_0^2 {\bf 1}_{\{|\xi_0| > \varepsilon
   n^{3/4} /  \max_j N_{n}(j)\}}\Big| {\cal G} \right) \frac{\alpha(n,0)}{n^{3/2}}=:\Sigma_{1}(\varepsilon,n) 
\end{array}$$
Let $\eta>0$.  We decompose the expectation of $\Sigma_1(\varepsilon,n)$ as the sum of
$$\Sigma_{1,1}(\varepsilon,n):= \mathbb{E}\Big( {\bf 1}_{\{n^{-3/4} \max_j N_{n}(j)\geq \eta\}} \mathbb{E} \left( \xi_0^2 {\bf 1}_{\{|\xi_0| > \varepsilon
   n^{3/4} /  \max_j N_{n}(j)\}}\Big| {\cal G} \right) \frac{\alpha(n,0)}{n^{3/2}}\Big)$$
   and
  $$\Sigma_{1,2}(\varepsilon,n):=\mathbb{E}\Big( {\bf 1}_{\{n^{-3/4}\max_j N_{n}(j) < \eta\}} \mathbb{E} \left( \xi_0^2 {\bf 1}_{\{|\xi_0| > \varepsilon
   n^{3/4} /  \max_j N_{n}(j)\}}\Big| {\cal G} \right) \frac{\alpha(n,0)}{n^{3/2}}\Big).$$ 
Using (ii)(a) of Proposition \ref{cadre}, $\Sigma_{1,2}(\varepsilon,n)$ is bounded by
$$C \mathbb{E} \left( \xi_0^2 {\bf 1}_{\{|\xi_0| > \varepsilon /\eta\}}\right)$$ 
From assumption $(A_1)$ of Theorem \ref{thmix0}, for any $\kappa>0$, there exists $\eta_0>0$ such that the above term is less than $\kappa/2$ for any $\eta\leq \eta_0$.\\*
We now fix $\eta$ equal to $\eta_0$. Using Cauchy-Schwarz inequality, 
\begin{eqnarray*} 
\Sigma_{1,1}(\varepsilon,n)&\leq & \mathbb{E}(\xi_0^2)\  \mathbb{P}
\Big(n^{-3/4} \max_j N_{n}(j)\geq \eta_0\Big)^{1/2}  \mathbb{E}\Big(\frac{\alpha(n,0)^2}{n^{3}}\Big)^{1/2}
\end{eqnarray*}
From assumption $(A_1)$ of Theorem \ref{thmix0}, (i) and (ii)(a) of Proposition \ref{cadre}, it follows that 
$\Sigma_{1,1}(\varepsilon,n)\leq \kappa/2$ for $n$ large enough.

Since $\varepsilon\rightarrow \mathbb{E}\big(\Sigma_{1}(\varepsilon,n)\big)$ is decreasing, we can find a sequence of positive numbers $(\varepsilon_n)_{n
  \geq 1}$ such that $\varepsilon_n \xrightarrow[n \rightarrow + \infty]{}
  0$, and 
\begin{equation}\label{thebigone}
  \mathbb{E}\Big(\Sigma_1(\varepsilon_n,n)\Big)\xrightarrow[n \rightarrow + \infty]{} 0
\end{equation}
Let us now prove that it implies
\begin{equation}\label{lindbisbis}
n^{-3/2} \bbE \left(\Big( \sum_{i\in\mathbb{Z}}  \varphi^{\varepsilon_n
  n^{3/4}}(X_{n,i})-\bbE(\varphi^{\varepsilon_n
  n^{3/4}}(X_{n,i})|{\cal G})\Big)^2 \right) \xrightarrow[n \rightarrow + \infty]{}
  0 \, .
\end{equation}

For any fixed $n \geq 0$, and any $i\in\mathbb{Z}$ such that $N_{n}(i)
\neq 0$, define~:
$$V_{n,i}=\varphi^{\varepsilon_n
    n^{3/4} / N_{n}(i)}(\xi_i) - \bbE(\varphi^{\varepsilon_n
    n^{3/4} / N_{n}(i)}(\xi_i)|{\cal G})$$
If $N_{n}(i)=0$, let $V_{n,i} =0$. As
    for any fixed $n \geq 0$ and any $i\in\mathbb{Z}$, the function
$$x \mapsto \varphi^{\varepsilon_n
    n^{3/4} / N_{n}(i)}(x) $$
is $1$-Lipschitz, we have for any fixed path of the random walk, for all $l \geq 1$, for all $k \geq 1$,
$$\theta_{k,2}^{V_{\cdot ,n}}(l) \leq \theta_{k,2}^{\xi}(l) \, ,$$
where $V_{\cdot,n}=(V_{n,i})_{i\in\mathbb{Z}}$ and
$\xi=(\xi_i)_{i \in
  \mathbb{Z}}$.
\begin{eqnarray*}
& &\mathbb{E}\Big(\Big(\sum_{j\in\mathbb{Z}} N_{n}(j) V_{n,j}\Big)^2\Big|{\cal G}\Big)\\
&=&\sum_{j\in\mathbb{Z}}  N_{n}(j)^2 \mathbb{E}(V_{n,j}^2|{\cal G}) +
\sum_{i\in\mathbb{Z}}\sum_{j\in\mathbb{Z}; j\neq i} N_{n}(i)\ N_{n}(j) \mathbb{E}(V_{n,i}V_{n,j}|{\cal G})\\
&\leq & \sum_{j\in\mathbb{Z}}  N_{n}(j)^2 \mathbb{E}(V_{n,j}^2|{\cal G}) + \sum_{i\in\mathbb{Z}}N_{n}(i)^2\sum_{j\in\mathbb{Z}; j\neq i}
|\mathbb{E}(V_{n,i}V_{n,j}|{\cal G})|
\end{eqnarray*}
by remarking that $N_{n}(i)\ N_{n}(j) \leq \frac{1}{2} (N_{n}(i)^2+N_{n}(j)^2)$.\\*
Then for any $j>i$, using Cauchy-Schwarz inequality, we obtain that
\begin{eqnarray*}
|\mathbb{E}(V_{n,i} V_{n,j}|{\cal G})|&=& \left|\mathbb{E}\left(V_{n,i} \; \mathbb{E} \left( V_{n,j}|{\cal M}_i\right)\Big|{\cal G}\right)\right| \\
&\leq & \mathbb{E}(V_{n,i}^2 | {\cal G})^{1/2} \ 
 \mathbb{E} \Big(\mathbb{E}(V_{n,j} | {\cal M}_i)^2 \Big|{\cal G}\Big)^{1/2}\\
&\leq & \mathbb{E}(V_{n,i}^2 | {\cal G})^{1/2} \ \theta_{1,2}^{\xi}(j-i).
\end{eqnarray*}
Moreover, as $N_{n}(i) V_{n,i}=\varphi^{\varepsilon_n
    n^{3/4} }(X_{n,i}) - \mathbb{E}(\varphi^{\varepsilon_n
    n^{3/4} }(X_{n,i})|{\cal G}),$ we get    
\begin{equation}\label{majbis}
n^{-3/2}\bbE\Big(\Big(\sum_{i\in\mathbb{Z}} \varphi^{\varepsilon_n
    n^{3/4} }(X_{n,i}) -\bbE(\varphi^{\varepsilon_n
    n^{3/4} }(X_{n,i})|{\cal G})\Big)^2\Big) \leq \mathbb{E}\Big(C_n \Big(\frac{\alpha(n,0)}{n^{3/2}}\Big)\Big)
\, ,
\end{equation}
with $C_n=\sup_{i\in\mathbb{Z} } \left(\mathbb{E}(V_{n,i}^2|{\cal G})\right)
+ 2
    \sqrt{\sup_{i\in\mathbb{Z}} \left(\mathbb{E}(V_{n,i}^2|{\cal G})\right)} \ \sum_{l=1}^{\infty}
\theta_{1,2}^{\xi}(l)$. 
It remains to prove that the right hand term in
    (\ref{majbis}) converges to $0$ as
    $n$ goes to infinity.

We have, for $n$ large enough,
\begin{equation}\label{maj2}
\mathbb{E}\Big(V_{n,i}^2\Big|{\cal G}\Big) \leq \mathbb{E} \left( \xi_i^2 {\bf 1}_{\{|\xi_i|
>  \varepsilon_n  n^{3/4} /  \max_j N_{n}(j)\}}\Big| {\cal G} \right) \, ,
\end{equation}
so using Cauchy-Schwarz inequality,
$$\mathbb{E}\Big(C_n\Big(\frac{\alpha(n,0)}{n^{3/2}}\Big)\Big)\leq  \mathbb{E}\Big(\Sigma_1(\varepsilon_n,n)\Big)+2\left(\sum_{l=1}^{\infty}
\theta_{1,2}^{\xi}(l)\right)\  \mathbb{E}\Big(\Sigma_1(\varepsilon_n,n)\Big)^{1/2} \mathbb{E}\Big(\frac{\alpha(n,0)}{n^{3/2}}\Big)^{1/2} ,
 $$
 which tends to 0 from (ii)(a) of Proposition \ref{cadre}, assumption $(A_2)$ from Theorem \ref{thmix0} and $(\ref{thebigone})$.

\vspace{1.5em}

Define $Z_{n,i}$ by
$$\varphi_{\varepsilon_n n^{3/4}}(X_{n,i})- \mathbb{E}\left(\varphi_{\varepsilon_n n^{3/4}}(X_{n,i})|{\cal G}\right) .$$
By (\ref{lindbisbis}) we conclude that, to prove Theorem \ref{thmix0}, it is
enough to prove it for the truncated sequence $\left(Z_{n,i}\right)_{n\geq
  0, i\in\mathbb{Z} }$, that is to show that
  \begin{equation}\label{tlctronc}
  \frac{1}{n^{3/4}}\sum_{i\in\mathbb{Z}}Z_{n,i}
\xrightarrow[ n \rightarrow +
  \infty]{\mathcal{D}} \sqrt{\sum_k \mathbb{E}(\xi_0 \xi_k)}\  \Delta_1 \, .
  \end{equation}

The proof is now a variation on the proof of Theorem 4.1 in
Utev (1990). Let
$$d_t(X,Y)=\left|\mathbb{E}(e^{itX}|{\cal G}) - \mathbb{E}(e^{itY}|{\cal G}) \right| \, .$$
Let $\eta$ be a random variable with standard normal distribution, independent of the random walk $(S_k)_{k\geq 0}$. Let $(X_i)_{i=1,..,3}$ be random
  variables such that for $i=1,..,3$, $\bbE(X_i|{\cal G})=0$ and $\bbE(X_i^2|{\cal G})$ is bounded. Let $Y_1,Y_2$ be random variables such that $\bbE(Y_i|{\cal G})=0$ and  $\bbE(Y_i^2|{\cal G})$ is bounded for $i=1,2$. They are assumed to be independent conditionnally to the random walk.
We define $$A_t(X)=d_t\left(X , \eta \sqrt{\mathbb{E}(X^2|{\cal G})}\right).$$
We first need some simple
  properties of $d_t$ and $A_t$~:
\begin{lem}[Lemma 4.3 in Utev, 1990]\label{ateate}
$$A_t(X_1)\leq \frac{2}{3}|t|^3 \mathbb{E}(|X_1|^3|{\cal G}),$$
$$A_t(Y_1+Y_2) \leq A_t(Y_1)+A_t(Y_2),$$
$$d_t(X_2+X_3,X_2) \leq \frac{t^2}{2} \left( \mathbb{E}(X_3^2|{\cal G}) +
  (\mathbb{E}(X_2^2|{\cal G}) \mathbb{E} (X_3^2|{\cal G}))^{1/2}\right),$$
$$d_t(\eta a, \eta b) \leq \frac{t^2}{2} |a^2-b^2|.$$
\end{lem}

We next need the following lemma~:
\begin{lem}\label{ut}
Let $0 < \varepsilon < 1$.
There exists some positive constant $C(\varepsilon)$ such that for all $a \in
\mathbb{Z}$, for all $v \in \mathbb{N}^*$, $A_t\left(n^{-3/4}\sum_{i=a+1}^{a+v}
Z_{n,i} \right)$ is bounded by
$$ C(\varepsilon)\left( |t|^3
  h^{2/\varepsilon} n^{-9/4} \sum_{i=a+1}^{a+v} \mathbb{E}\left(|Z_{n,i}|^3|{\cal G}\right) +
  t^2 \left(h^{\frac{\varepsilon -1}{2}} + \sum_{j~: 2^j \geq h^{1/\varepsilon}}
  2^{\frac{3 j}{2}}g(2^{j\varepsilon})\right) n^{-3/2}
  \sum_{i=a+1}^{a+v}N_{n}(i)^2\right),$$
where $h$ is an arbitrary positive natural number and with $g$ introduced in
  Assumption $(A_2)$ of Theorem \ref{thmix0}.
\end{lem}

\vspace{1.5em}

Before proving Lemma \ref{ut}, we achieve the proof of Theorem \ref{thmix0}.
We can decompose 
\begin{equation}\label{in}
\mathbb{E}\left(e^{it\frac{1}{n^{3/4}}\sum_{i\in\mathbb{Z}} Z_{n,i}}\right) - \mathbb{E}\left(e^{it \sqrt{\sum_k \mathbb{E}(\xi_0 \xi_k)} \Delta_1}\right)
\end{equation}
as the sum of $I_i(n), i=1,..,4$ where 
\begin{eqnarray*}
I_1(n)&= & \mathbb{E}\left(e^{it\frac{1}{n^{3/4}}\sum_{i\in\mathbb{Z}} Z_{n,i}}\right) -\mathbb{E}\left(e^{-\frac{t^2}{2n^{3/2}} \mathbb{E}\left(\left(\sum_{i\in\mathbb{Z}}
Z_{n,i}\right)^2|{\cal G}\right) }\right) \\
I_2(n) &=& \mathbb{E}\left(e^{-\frac{t^2}{2n^{3/2}} \mathbb{E}\left(\left(\sum_{i\in\mathbb{Z}}
Z_{n,i}\right)^2|{\cal G}\right) }\right) - 
 \mathbb{E}\left(e^{-\frac{t^2}{2n^{3/2}} \mathbb{E}\left(\left(\sum_{i\in\mathbb{Z}}
X_{n,i}\right)^2|{\cal G}\right) }\right)\\
I_3(n)&=&  \mathbb{E} \left(e^{-\frac{t^2}{2n^{3/2}} \mathbb{E}\left(\left(\sum_{i\in\mathbb{Z}}
X_{n,i}\right)^2|{\cal G}\right) }\right)-
 \mathbb{E}\left(e^{-\frac{t^2}{2 n^{3/2}} \sum_{i,j\in\bbZ} N_n(i)^2\mathbb{E}(\xi_i \xi_j)} \right)\\
I_4(n)&=&\mathbb{E}\left(e^{-\frac{t^2}{2 n^{3/2}} \sum_{i,j\in\bbZ} N_n(i)^2\mathbb{E}(\xi_i \xi_j)}\right) - \mathbb{E}\left(e^{it \sqrt{\sum_k \mathbb{E}(\xi_0 \xi_k)} \Delta_1}\right)
\end{eqnarray*}
To prove Theorem \ref{thmix0}, it is enough to prove that for any $i=1,..,4$, $I_i(n)$ goes to 0 as $n\rightarrow +\infty$.\\*
\noindent{\bf Estimation of $I_1(n)$~:} Let us denote by $M_n$ the random variable $\max_{k=0,\ldots,n} |S_k|$. From Lemma \ref{ut}, we have 
$$|I_1(n)| \leq \mathbb{E} \left(A_t\left(n^{-3/4}\sum_{i=-M_n}^{M_n}
Z_{n,i}\right)\right)
 \leq  C(t,\varepsilon) \left( h^{2 / \varepsilon }  n^{-9/4} \sum_{i\in\bbZ} \bbE(|Z_{n,i}|^3)
+ \delta(h)\right) ,$$ 
with $\delta(h)=\Big(h^{\frac{\varepsilon
-1}{2}}+\sum_{j~: 2^j \geq h^{\frac{1}{\varepsilon}}}
  2^{\frac{3 j}{2}}g(2^{j\varepsilon})\Big) 
  \mathbb{E}\left(\frac{\alpha(n,0)}{n^{3/2}}\right)$.

\vspace{1em}

Hence using assumption $(A_2)$ from Theorem \ref{thmix0} and
 $(ii)(a)$ from Proposition \ref{cadre}, we get $\delta(h)\xrightarrow[ h
   \rightarrow + \infty]{} 0 $.

On the other hand, from assumption $(A_1)$ of Theorem \ref{thmix0}, there exists a constant $C>0$ such that
\begin{equation}\label{hip}
n^{-9/4}\sum_{i\in\bbZ} \bbE(|Z_{n,i}|^3) \leq C \varepsilon_n \mathbb{E}\left(\frac{\alpha(n,0)}{n^{3/2}}\right)
\end{equation}
 which tends to zero as $n$ tends to
    infinity, using $(ii)(a)$ from Proposition \ref{cadre} and the fact that $\varepsilon_n
    \xrightarrow[n \rightarrow + \infty]{}0$.
Consequently
$$\inf_{h \geq 1} \left(h^{2 / \varepsilon }\sum_{i\in\mathbb{Z}} \mathbb{E}(|Z_{n,i}|^3)
+ \delta(h)\right) \xrightarrow[n \rightarrow + \infty]{}0.$$

\noindent{\bf Estimation of $I_2(n)$~:} Using that for any $x,y\geq 0$, $|e^{-x}-e^{-y}|\leq |x-y|$ and the fact that 
$$Z_{n,i}=X_{n,i}-[\varphi^{\varepsilon_n n^{3/4}}(X_{n,i})-\mathbb{E}(\varphi^{\varepsilon_n n^{3/4}}(X_{n,i})|{\cal G})],$$
 we deduce that
\begin{eqnarray*}
|I_2(n)|&\leq &\frac{t^2}{n^{3/2}}\left[\mathbb{E} \left(\Big(\sum_{i\in\bbZ} \varphi^{\varepsilon_n n^{3/4}}(X_{n,i})-\mathbb{E}(\varphi^{\varepsilon_n n^{3/4}}(X_{n,i})|{\cal G}) \Big)^2\right)\right.\\
&+&\left.\mathbb{E}\left(\Big(\sum_{i\in\bbZ} X_{n,i}\Big)^2\right)^{1/2} 
\mathbb{E} \left(\Big(\sum_{i\in\bbZ} \varphi^{\varepsilon_n n^{3/4}}(X_{n,i})-\mathbb{E}(\varphi^{\varepsilon_n n^{3/4}}(X_{n,i})|{\cal G}) \Big)^2\right)^{1/2}\right]
\end{eqnarray*}
Since we have
\begin{eqnarray*}
\mathbb{E}\Big(\Big(\sum_{i\in\bbZ}
X_{n,i}\Big)^2\Big)&=& \sum_{i,j\in \mathbb{Z}} \mathbb{E}(N_n(i) N_n(j)) \mathbb{E}(\xi_i \xi_j)\\
&=& \sum_{i\in \mathbb{Z}} \mathbb{E}(\alpha(n,i)-\alpha(n,0)) \mathbb{E}(\xi_0 \xi_i) +  \mathbb{E}(\alpha(n,0)) \sum_{i\in\mathbb{Z}} \mathbb{E}(\xi_0 \xi_i)\\
&\leq & C n^{3/2}  
\end{eqnarray*}
by combining (ii)(a)--(iii) of Proposition \ref{cadre} and $(\ref{imp})$. 
Then, from (\ref{lindbisbis}), we deduce that $I_2(n)$
converges to 0.

\noindent{\bf Estimation of $I_3(n)$~:} 
We have 
$$|I_3(n)|\leq \frac{t^2}{2n^{3/2}} \sum_{i\in\mathbb{Z}} \mathbb{E}(|\alpha(n,i)-\alpha(n,0)|) |\mathbb{E}(\xi_0 \xi_i)| =o(1),$$
by combining (iii) of Proposition \ref{cadre} and $(\ref{imp})$.

\noindent{\bf Estimation of $I_4(n)$~:} 
From (ii)(b) of Proposition \ref{cadre}, we know that the sequence  $$\mathbb{E}\left(e^{-\frac{t^2}{2 n^{3/2}} \sum_{i,j\in\bbZ} N_n(i)^2\mathbb{E}(\xi_i \xi_j)}\right)$$
 converges to 
$ \mathbb{E}\left(e^{-\frac{t^2}{2} \sum_{i\in\bbZ} \mathbb{E}(\xi_0 \xi_i) \int_{\bbR} L_1^{2}(x) \, dx} \right)$
which is equal to the characteristic function of the random variable $\sqrt{\sum_{i\in\bbZ} \mathbb{E}(\xi_0 \xi_i)} \Delta_1.$

\vspace{1.5em}

\noindent {\bf Proof of Lemma \ref{ut}:} Let $h \in \mathbb{N}^*$.
Let $0 < \varepsilon < 1$.  In the following, $C$,
$C(\varepsilon)$ denote constants which may vary from line to line.
Let $\kappa_{\varepsilon}$ be a positive constant greater than $1$
which will be precised further. Let $v < \kappa_{\varepsilon}\
h^{\frac{1}{\varepsilon}}$. We have
\begin{eqnarray}\label{petit}
A_t \left(n^{-3/4}\sum_{i=a+1}^{a+v} Z_{n,i}\right)  &\leq &\frac{2|t|^3}{3n^{9/4}}
\mathbb{E}\left(\left|\sum_{i=a+1}^{a+v} Z_{n,i}\right|^3\Big|{\cal G}\right) \nonumber \\*
&\leq &\frac{2}{3n^{9/4}}  \kappa_{\varepsilon}^2 \ |t|^3 \
h^{2/\varepsilon}\sum_{i=a+1}^{a+v} \mathbb{E}(|Z_{n,i}|^3|{\cal G})
\end{eqnarray}
since $|x|^3$ is a convex function.

Let now $v \geq \kappa_{\varepsilon}\
h^{\frac{1}{\varepsilon}}$. Without loss of generality, assume that $a=0$. Let 
$\delta_{\varepsilon}=(1-\varepsilon^2 +2\varepsilon)/2$. Define then
$$m=[v^{\varepsilon}] , \ B= \left\{u \in \mathbb{N}~: \ 2^{-1} (v
-[v^{\delta_{\varepsilon}}]) \leq um \leq 2^{-1}v\right\},$$
$$A=\left\{ u \in \mathbb{N}~:0 \leq u \leq v, \
\sum_{i=um+1}^{(u+1)m}N_{n}(i)^2 \leq (m/v)^{\varepsilon} \sum_{i=1}^v
N_{n}(i)^2 \right\}.$$
Following Utev (1991) we prove that, for $0 < \varepsilon < 1$, $A \cap B$
is not wide for $v$ greater than $\kappa_{\varepsilon}$. We have indeed
$$|A \cap B| =|B|-|\overline{A}\cap B|\geq |B|-|\overline{A}|\geq
\frac{v^{(1-\varepsilon^2)/2}}{2}\Big(1-4v^{-(1-\varepsilon)^2/2}\Big)-\frac{3}{2} \ , $$
where $\overline{A}$ denotes the complementary of the set $A$.
We can find $\kappa_{\varepsilon}$ large enough so that $|A \cap B|$ be positive.

Let $u \in A \cap B$. We start from
the following simple identity
\begin{eqnarray}
\nonumber Q & \equiv & n^{-3/4}\sum_{i=1}^v
Z_{n,i}\\
\nonumber & = &
n^{-3/4}\sum_{i=1}^{um}Z_{n,i}+n^{-3/4}\sum_{i=um+1}^{(u+1)m}Z_{n,i}+n^{-3/4}\sum_{i=(u+1)m+1}^v
Z_{n,i}\\
\label{4.7} & \equiv & Q_1 + Q_2 + Q_3.\end{eqnarray}For any fixed $n \geq 0$, and any $i \in\mathbb{Z}$ such that $N_{n}(i)
\neq 0$, define~:
$$W_{n,i}=\varphi_{\varepsilon_n n^{3/4} / N_{n}(i)}(\xi_i)-\mathbb{E}\left( \varphi_{\varepsilon_n n^{3/4} / N_{n}(i)}(\xi_i)\Big| {\cal G} \right) \, .$$
If $N_{n}(i)=0$, let $W_{n,i} =0$. As
    for any fixed $n \geq 0$, $i\in\mathbb{Z}$, the function
$$x \mapsto \varphi_{\varepsilon_n
    \sigma_n / N_{n}(i)}(x)-\mathbb{E}\left( \varphi_{\varepsilon_n
    \sigma_n / N_{n}(i)}(\xi_i)\Big| {\cal G}  \right)$$
is $1$-Lipschitz, we have for any fixed path of the random walk, for all $l \geq 1$, for all $k \geq 1$,
$$\theta_{k,2}^{W_{\cdot ,n}}(l) \leq \theta_{k,2}^{\xi}(l) \, ,$$
where $W_{\cdot,n}=(W_{n,i})_{i\in\mathbb{Z}}$ and
$\xi=(\xi_i)_{i \in
  \mathbb{Z}}$.
Hence, arguing as for the proof of inequality (\ref{majbis}), we prove that for any fixed $n$, any $a,b\in\bbN$,
\begin{equation}\label{maj3}
\bbE\Big(|\sum_{i=a}^{b}Z_{n,i}|^2\ \Big| \, {\cal G}\Big)\leq C \sum_{i=a}^{b} N_{n}(i)^2
\end{equation}
with $C =2\mathbb{E}(\xi_0^2) + 2 \sqrt{2} \mathbb{E}(\xi_0^2)^{1/2} \ \sum_{l=1}^{\infty}
\theta_{1,2}(l)$ which is finite from assumptions $(A_1)$ and
$(A_2)$.
By Lemma \ref{ateate},
\begin{equation}\label{4.8}
d_t(Q,Q_1+Q_3)=d_t(Q,Q-Q_2)\leq \frac{t^2}{2} \left(\mathbb{E}( Q_2^2|{\cal G}) +
\mathbb{E}(Q_2^2|{\cal G})^{1/2} \mathbb{E}(Q^2|{\cal G})^{1/2}\right).
\end{equation}
Using (\ref{4.8}) and (\ref{maj3}), we get
\begin{equation}\label{4.9}
d_t(Q,Q_1+Q_3) \leq C t^2 v^{\frac{(\varepsilon-1)\varepsilon}{2}} n^{-3/2} \sum_{i=1}^v N_{n}(i)^2.
\end{equation}
Now, given the random variables $Q_1$ and $Q_3$, we define two random variables $g_1$ and $g_3$ which are assumed independent conditionally to the random walk $(S_k)_{k\geq 0}$ such that conditionnally to the random walk, the distribution of $g_i$ coincides with that of $Q_i$, $i=1,3$.
We have
$$\begin{array}{rcl}
& & d_t(Q_1+Q_3,g_1+g_3)\\*
 & =&\left|\mathbb{E}((e^{itQ_1}-1)(e^{itQ_3}-1)|{\cal G}) -\mathbb{E}(e^{itQ_1}-1|{\cal G})\mathbb{E}(e^{itQ_3}-1|{\cal G})\right|\\\\
& \leq  &  \mathbb{E}(|e^{itQ_1}-1|^2|{\cal G})^{1/2}
\mathbb{E}\big(\big|\mathbb{E}\left(e^{itQ_3}-1-\mathbb{E}(e^{itQ_3}-1|{\cal G}) \ | \,
\mathcal{M}_{um};{\cal G}\right)\big|^{2}|{\cal G}\big) \\\\
& \leq & 2 |t| n^{-3/4}\bbE\Big(|\sum_{i=1}^{um}Z_{n,i}|^2\ | \, {\cal G}\Big)^{1/2} v \ |t| \ n^{-3/4}\left(
        \sum_{i=(u+1)m+1}^v
N_{n}(i)\right) \theta_2^{\xi}(m+1)\\
&\leq & C \ t^2 \ v^{3/2} \ n^{-3/2} \left(\sum_{i=1}^v
N_{n}(i)^2\right) g(v^{\varepsilon})
\end{array}$$
by (\ref{maj3}), Definition \ref{coeffsuite} and Assumption $(A_2)$ of
Theorem \ref{thmix0}.
Hence
\begin{equation}\label{4.10}
d_t(Q_1+Q_3,g_1+g_3) \leq  C t^2 f(v) n^{-3/2}
\sum_{i=1}^v N_{n}(i)^2,
\end{equation}
where $f(v)=v^{3/2} \ g(v^{\varepsilon})$ is non-increasing by assumption
$(A_2)$ of Theorem \ref{thmix0}.

We also have by Lemma \ref{ateate}
\begin{equation}\label{4.11}
A_t(g_1+g_3)\leq A_t(g_1)+A_t(g_3).
\end{equation}
Finally, still by Lemma \ref{ateate}, and using Definition \ref{coeffsuite}, we have
\begin{eqnarray}
& &\nonumber d_t \left( \eta \sqrt{\mathbb{E}(Q^2|{\cal G})}, \eta
\sqrt{\mathbb{E}\left( (g_1+g_3)^2|{\cal G} \right) } \right)\\
 & \leq & \frac{t^2}{2}
\left| \mathbb{E} (Q^2|{\cal G}) - \mathbb{E} \left( (g_1+g_3)^2|{\cal G} \right)  \right|\\
 \nonumber & \leq & \frac{t^2}{2} \left| \mathbb{E} (Q_2^2|{\cal G}) + 2 \mathbb{E}
(Q_1Q_2|{\cal G}) + 2 \mathbb{E}(Q_2Q_3|{\cal G}) +2 \mathbb{E}(Q_1Q_3|{\cal G})\right|\\
\label{4.12} & \leq & C t^2 \left(  v^{\frac{(\varepsilon -1)\varepsilon}{2}} +f(v)\right) n^{-3/2}\sum_{i=1}^v N_{n}(i)^2.
\end{eqnarray}
Combining (\ref{4.9})-(\ref{4.12}), we get the following recurrent
inequality~:
$$\begin{array}{rcl}
A_t\left(n^{-3/4}\sum_{i=1}^v Z_{n,i}\right) & \leq & A_t\left(n^{-3/4}\sum_{i=1}^{um}
Z_{n,i}\right) + A_t\left(n^{-3/4}\sum_{i=(u+1)m+1}^{v} Z_{n,i}\right)\\
& + & C t^2 \left(  v^{\frac{(\varepsilon -1)\varepsilon}{2}} + f(v)\right) n^{-3/2}
\sum_{i=1}^v N_{n}(i)^2
\end{array}$$
for $v \geq \kappa_{\varepsilon} \
h^{\frac{1}{\varepsilon}}\geq \kappa_{\varepsilon}$.

We then need the following Lemma, which is a variation on Lemma 1.2. in Utev
(1991).

\begin{lem}\label{ututut}
For every $\varepsilon\in\ ]0,1[$, denote $\delta_{\varepsilon}=(1-\varepsilon^2+2\varepsilon)/2$.
Let a non-decreasing sequence of non-negative numbers $a(n)$ be
specified, such that there exist non-increasing sequences of
non-negative numbers $\varepsilon(k)$, $\gamma(k)$ and a sequence of
naturals $T(k)$, satisfying conditions
$$T(k) \leq 2^{-1}(k+[k^{\delta_{\varepsilon}}]),$$
$$a(k) \leq \max_{k_0 \leq s \leq k}\left(a(T(s))+\gamma(s)\right)$$
for any $k \geq k_0$ with an arbitrary $k_0 \in \mathbb{N}^*$. Then
$$a(n) \leq a(n_0)+2 \sum_{k_0 \leq 2^j \leq n}\gamma(2^j) \ ,$$
for any $n \geq k_0$, where one can take $n_0=2^c$ with $c>\frac{2-\delta_{\varepsilon}}{1-\delta_{\varepsilon}}$.
\end{lem}

\vspace{1.5em}

\noindent {\bf Proof of Lemma \ref{ututut}:} The proof follows essentially
the same lines as the proof of Lemma 1.2. in Utev (1991) and therefore is
omitted here. \epreuve

\vspace{1.5em}

We now apply Lemma \ref{ututut} above with
\begin{itemize}
\item[$\star$] $k_0=\kappa_{\varepsilon} \ h^{\frac{1}{\varepsilon}}$,
\item[$\star$] for $k \ge k_0$, $T(k)=\max \left\{ u_km_k, k-u_km_k-m_k \right\}$ where $u_k$
  and $m_k$ are defined from $k$ as $u$ and $m$ from $v$ (see the proof of
  $A \cap B $ not wide),
\item[$\star$] $c<\frac{\ln(\kappa_{\varepsilon})}{\ln(2)}\ \ \ $ (we may
  need to enlarge $\kappa_{\varepsilon}$),
\item[$\star$] for $s \geq k_0$, $\gamma(s)=C  \ t^2 \ \left( s^{\frac{\varepsilon ( \varepsilon
-1)}{2}}+f(s) \right)$,
\item[$\star$] for $s \geq k_0$, $a(s)=\displaystyle\sup_{l\in\mathbb{Z}} \ \max_{k_0 \leq i \leq s} \
\frac{A_t\left(n^{-3/4}\sum_{j=l+1}^{l+i}Z_{n,j}\right)}{n^{-3/2}\sum_{j=l+1}^{l+i}N_{n}(j)^2}
$.
\end{itemize}
Applying Lemma \ref{ututut} yields the statement of Lemma \ref{ut}. 
\\*
\\*
\noindent{\bf Proof of the tightness:}\\*
By Theorem 13.5 of Billingsley (1999), it suffices to prove
that there exists $K>0$ such that for all $t,t_1,t_2\in[0,T],
T<\infty,$ s.t. $t_{1}\leq t\leq t_{2},$ for all $n\geq 1$,
\begin{equation}\label{pro}
\bbE\Big(|\Sigma_{[nt]}-\Sigma_{[nt_1]}| \cdot \ |\Sigma_{[nt_2]}-\Sigma_{[nt]}|\Big)\leq K
n^{3/2} |t_{2}-t_{1}|^{\frac{3}{2}}.
\end{equation}
Using Cauchy-Schwarz inequality, it is enough to prove that there
exists $K>0$ such that for all $t_1\leq t$, for all $n\geq 1$,
\begin{equation}\label{pro1}
\bbE\Big(|\Sigma_{[nt]}-\Sigma_{[nt_1]}|^2\Big)\leq K
n^{3/2} |t-t_{1}|^{\frac{3}{2}}.
\end{equation}
Since $ab\leq (a^2+b^2)/2$ for any real $a,b$, we have 
\begin{eqnarray*}
\bbE\Big(|\Sigma_{[nt]}-\Sigma_{[nt_1]}|^2\Big)&\leq &\sum_{i\in\bbZ}|\bbE(\xi_0\xi_i)|\ \sum_{j\in\bbZ} \bbE\Big( (N_{[nt]}(j)-N_{[nt_1]}(j))^2\Big)\\
&=& C \bbE\Big( \alpha\big([nt]-[nt_1]-1,0\big) \Big)\leq C n^{3/2}|t-t_1|^{3/2}
\end{eqnarray*}
using (\ref{imp}), the Markov property for the random walk $S$ and (ii)-(a) of Proposition \ref{cadre}. 
\epreuve

\end{document}